\documentclass[11pt]{amsart}
\usepackage{amsfonts,amsmath,amssymb}
\usepackage{graphicx}
\usepackage{epsfig}

\numberwithin{equation}{section}
\newcommand{\eps}{\varepsilon}

\newcommand{\E}{\mathsf{E}}
\newtheorem{lemma}{Lemma}[section]

\newtheorem{theorem}{Theorem}[section]

\newcommand{\Pp}{\mathsf{P}}

\newcommand{\R}{\mathbb{R}}
\newcommand{\N}{\mathbb{N}}

\newcommand{\sgn}{\mathop{\rm sgn}}
\newcommand{\Law}{\mathop{\rm Law}}
\newcommand{\Nc}{\mathcal{N}}
\newcommand{\ONE}{{\bf 1}}

\newcommand{\bpf}[1][Proof]{{\noindent {\sc #1: }}}
\newcommand{\epf}{{{\hspace{4 ex} $\Box$ \smallskip}}}
\author{
  Yuri Bakhtin
}
\address{School of Mathematics, Georgia Tech, Atlanta GA, 30332-0160, USA}
\email{bakhtin@math.gatech.edu}
\title{Decision making times in mean-field dynamic Ising model}

\begin{document}
\begin{abstract} We consider a dynamic mean-field ferromagnetic model in the low-temperature
regime in the neighborhood of the zero magnetization state. We study the random time it takes 
for the system to make a decision, i.e., to exit the neighborhood of the unstable 
equilibrium and approach one of the two stable equilibrium points.
We prove a limit theorem for the distribution of this 
random time in the thermodynamic limit.
\end{abstract}

\maketitle

\section{Introduction}

The asymptotic properties of exit from a small neighborhood of an unstable equilibrium of
a dynamical system under small white noise perturbation were first studied rigorously in~\cite{Kifer}.
It was shown that as the noise intensity~$\eps$ goes to zero, the exit time $\tau_\eps$ behaves
roughly as $a^{-1}\ln \eps^{-1}$, where $a$ is the local expansion rate (Lyapunov exponent)
at the equilibrium point. In \cite{Bakhtin-SPA:MR2411523} and~\cite{nhn} it was shown that
$\tau_\eps-a^{-1}\ln \eps^{-1}$ converges to a limiting distribution that is a dilation and
translation of $\ln|G|^{-1}$, where $G$ is a standard Gaussian random variable. 

Understanding distributional asymptotics for the exit time was pivotal in describing the vanishing noise
asymptotics for noisy heteroclinic networks, see~\cite{nhn} and~\cite{nhn-ds}. These systems occur naturally
in the context of neural dynamics and sequential decision making, see, e.g.,~\cite{rabinovich-et-al-2008}
and references therein. Exit times for diffusion models have been used in psychology to describe reaction times
in decision tasks, see~\cite{ratcliff-2008} and references therein, 
and it is natural to ask if the limiting behavior of exit times described above 
is reproduced in statistical mechanics models of neural computation. 

In this paper we study one of the the simplest possible models of this kind, the dynamic
mean field ferromagnetic model, also known as the Curie-Weiss model, in the low temperature regime with
two minima of free energy. We start the evolution of the system at the completely disordered state with zero magnetization,
where the numbers of plus spins equals the number of minus spins.  We stop the dynamics as soon as magnetization
enters a neighborhood of one of the stable equilibrium values and interpret that event as a decision made by the system between
the two choices. We show that as the number $N$ of spin
variables  (representing individual neurons in the  neural computation context) goes to infinity, the exit time behaves as $\ln N$ and the correction
to the main term converges to an affine transformation of $\ln|G|^{-1}$, thus reproducing the above result for the 
diffusion in the neighborhood of an unstable equilibrium.

\section{The model and the main result}
First let us recall the  mean-field ferromagnetic equilibrium Ising model  also known as Curie--Weiss model, see
~\cite[Section IV.4]{Ellis:MR2189669}.
Let us fix a large number $N$ and consider $N$ spin variables. Each variable $X_k, k=1,\ldots,N$
takes values $\pm 1$, and the energy assigned to a configuration $(x_k)_{k=1}^N$ is given by
\[
 E(x)=-\frac{1}{2N}\sum_{i,j=1}^N x_i x_j.
\]

We then can fix an inverse temperature value $\beta>0$ and consider the Boltzmann--Gibbs distribution
defined by $E$ and $\beta$:
\[
 \Pp_N\{X_k=x_k,\quad k=1,\ldots N\}=\frac{e^{-\beta E(x)}}{Z_N},
\]
where
\[
 Z_N=\sum_{x\in \{-1,1\}^N} e^{-\beta E(x)}
\]
is the partition function.

Since there is no geometry involved in this mean-field model and
the strength of interactions between two spins is the
same for all pairs of spins, one can describe the macroscopic behavior of the system by a single
variable called magnetization,
\[
 M(x)=\frac{1}{N}\sum_{i} x_i\ \in[-1,1].
\]
Notice that
\[
 E(x)=-N M^2(x)/2
\]
Therefore,
\[
 \Pp_N\left\{M(X)=\frac{n}{N}\right\}=\frac{1}{Z_N}{\binom{N}{(N+n)/2}}e^{N\frac{\beta}{2}
\cdot\left(\frac{n}{N}\right)^2}
\]
if $(N+n)/2$ is integer.

Recall that (see, e.g., \cite[Lemma I.3.2]{Ellis:MR2189669}) uniformly in $k=0,\ldots,n$,
\[
 \frac{1}{N}\ln\binom{N}{k}=h\left(\frac{k}{N}\right)+O\left(\frac{\ln N}{N}\right),
\]
where
\[
 h(x)=-x\ln x-(1-x)\ln(1-x),\quad x\in[0,1],
\]
is the entropy of the Bernoulli distribution with probabilities $x$ and $1-x$.
Therefore,
\[
 \Pp_N\left\{M(X)=\frac{n}{N}\right\}= \frac{1}{Z_N} e^{-N F(n/N)+O(\ln N)},
\]
where the free energy per spin $F$ is defined by
\[
F(m)=-\frac{\beta}{2}m^2-h\left(\frac{1}{2}+\frac{m}{2}\right).
\]
It is easy to show that the sequence of distributions $\Pp_N\{M(X)\in \cdot\}$ satisfies a
large deviation principle on $[0,1]$ with rate function $J$ given by
\[
 J(m)=F(m)-\min_{[0,1]} F.
\]
Differentiating $F$, we see that the minimizers of $F$ satisfy
\begin{equation}
\label{eq:equation_on_critical_points}
 \beta m=\frac{1}{2}\ln\frac{1+m}{1-m},
\end{equation}
and, as elementary analysis shows, (i) for $\beta<1$,  a unique minimizer of $F$ is $m=0$ (corresponding to completely disordered
case), and $F''(0)>0$ so that $F$ is approximately quadratic in the
neighborhood of the minimizer;  (ii) for $\beta>1$, there are two minimizers $m=\pm m_*$, for some $m^*>0$; (iii) if $\beta=1$ then 0 is still a unique minimizer, but contrary to the first case, $F''(0)=0$, and the leading term in the Taylor
expansion of $F$ at 0 is order 4.

\medskip

In this paper we are concerned with the low-temperature case (ii).  In that situation, point 0 is also
a solution of~\eqref{eq:equation_on_critical_points}, but it is an unstable equilibrium of the system being the local maximum of the free energy $F$. We are going to consider stochastic dynamics compatible with
Curie--Weiss model and study it in the neighborhood of the unstable equilibrium in the case $\beta>1$.

We must study a $\{-1,+1\}^N$-valued Markov process with intensities of spin flips
$c_i(x), i\in\{1,\ldots,N\}, x\in\{-1,+1\}^N$ defined by
\[
 \Pp\{X_i(t+\Delta t)\ne X_i(t)|X(t)=x\}=c_i(x)\Delta t+o(\Delta t),\quad \Delta t\downarrow 0.
\]
If we want the process $X$ to be reversible w.r.t.\ the Gibbs distribution $P_N$, we have to require that
\[
 c_i(x)\exp\left\{-\beta E(x)\right\}
\]
does not depend on $x_i$, see \cite[Section IV.2]{Liggett:MR776231}. Equivalently, we can require that
\[
 c_i(x)\exp\left\{\beta x_i M(x)\right\}
\]
does not depend on $x_i$. There are many choices for rates $c_i$, and there is no physical reason to prefer one
of them to others. In this paper we will work with
\begin{equation}
\label{eq:determining_the_rates}
 c_i(x)=\exp\left\{-\beta x_i M(x)\right\},
\end{equation}
although our results should hold for a variety of other choices of $c_i(x).$ 

Notice that if the spin $x_i$ is aligned with magnetization $M(x)$, then the resulting flipping rate of $i$-th spin
is lower than in the opposite situation where $x_i$ is misaligned with $M(x)$. This is the result of the ferromagnetic
nature of the model which favors configurations with most spins aligned with each other.

Suppose now that we observe only the magnetization, or, equivalently, the number of $+1$-spins. Flipping a $-1$ spin means then
a transition from the current magnetization $m$ to $m+2/N$. Since the number of $-1$ spins equals $N(1-m)/2$, we see that the total
transition rate $m\mapsto m+2/N$ is given by $\lambda_+(m,N)$, where
\[
 \lambda_+(m,N)=N\frac{1-m}{2}\exp\{\beta m\}.
\]
Flipping a $+1$ spin means a transition from the current magnetization $m$ to $m-2/N$. Since the number of $+1$ spins equals $N(1+m)/2$, we see that the total
transition rate $m\mapsto m-2/N$ is given by $\lambda_-(m,N)$, where
\[
 \lambda_-(m,N)=N\frac{1+m}{2}\exp\{-\beta m\}.
\]

Let us consider the Markov process $M_N$ describing the evolution of magnetization in the above model and set $M_N(0)=0$
(this means that $N$ has to be even, but this is not a really important restriction).

It is clear that $M_N$ will spend some time in the neighborhood of $0$ and then it will escape that neighborhood and head
towards one of the minims of free energy, $\pm m_*$. We can interpret the exit in each of these directions as the decision
made by the system. We set a threshold level $R\in(0,m_*)$ and as soon as $M_N$ exceeds $R$ in absolute value, we claim that
the system has made the decision. The choice of one of the two alternatives is encoded by the sign of $M_N$ at that time.
Our main result describes the asymptotics of the random time it takes to reach the threshold $R$ starting
from the completely disordered state with zero magnetization. According to the interpretation above, this time can be viewed 
as the decision making time for the situation where the initial state is a completely unbiased indecisive state. 

More formally,  for any $R\in(0,m_*)$ we introduce
\[
 \tau_N(R)=\inf\{t: |M_N(t)|\ge R\}.
\]
Our main result describes the joint asymptotic behavior of random variables~$\tau_N(R)$ and $\sgn M_N(\tau_N)$. To state it, we need more notation.
For $m\in[-1,1]$, we denote 
\begin{equation}
 \label{eq:drift}
 b(m)=\frac{2}{N}(\lambda_+(m,N)-\lambda_-(m,N)) =(1-m)e^{\beta m} - (1+m)e^{-\beta m},
\end{equation}
introduce
$a=b'(0)=2\beta-2>0$ and  $Q(x)=b(x)-ax, x\in\R$,
and define
\[
 D(R)=K(R)+\frac{\ln R}{a}+\frac{\ln (a/2)}{2a},\quad R\in\R, 
\]
where
\begin{equation}
\label{eq:K(R)} 
K(R)=-\int_0^{R}\frac{Q(x)}{axb(x)}dx\in\R,\quad R\in\R.
\end{equation}

\begin{theorem} 
 For any $R\in(0,m^*)$, as $N\to\infty$
\[
 \left(\sgn M_N(\tau_N(R)),\ \tau_N(R) -\frac{1}{2a}\ln N\right)\ \stackrel{distr}{\longrightarrow}\ \left (\sgn G, -\frac{1}{a}\ln|G|+D(R)\right),
\]
where $G$ is a standard Gaussian random variable.
\end{theorem}

\section{Proof}
The proof is based on the theory of Markov processes, martingales and their convergence. We refer to
\cite{Ethier-Kurtz:MR838085} as an excellent source on the relevant background.

The form of rates $\lambda_+,\lambda_-$ implies that 
the generator $L_N$ of the magnetization process $M_N$ is given by
\begin{align*}
 L_Nf(m)=\frac{N}{2}&\left[(1-m)e^{\beta m}\left(f(m+\frac{2}{N})-f(m)\right)\right.\\
&\left.+(1+m)e^{-\beta m}\left(f(m-\frac{2}{N})-f(m)\right)\right].
\end{align*}
In particular,
\[
 f(M_N(t))-\int_0^tL_Nf(M_N(s))ds
\]
has to be a martingale for any $f$, see Proposition~1.7 in~\cite[Chapter 4]{Ethier-Kurtz:MR838085}. Choosing $f(x)\equiv x$ on $[-1,1]$, we obtain
that
\[
 Z_N(t)=M_N(t)-\int_0^tb(M_N(s))ds
\]
is a bounded variation cadlag martingale, where $b(\cdot)$ is defined in~\eqref{eq:drift}, so that it plays the role of drift coefficient
that drives the deterministic component of the process. Notice that zeros of $b(m)$ coincide with solutions of equation~\eqref{eq:equation_on_critical_points},
so that for any point $x\in (0,m_*)$, $\lim_{t\to+\infty} S^tx =m_*$ and $\lim_{t\to-\infty} S^tx =0$, and for any point $x\in (-m_*,0)$, 
$\lim_{t\to+\infty} S^tx =-m_*$ and $\lim_{t\to-\infty} S^tx =0$, where $S$ is the flow generated by $b$:
$S^0x=x$ and $\frac{d}{dt} S^tx=b(S^tx)$.

Using the representation
\begin{equation}
\label{eq:Q}
 b(m)=am + Q(m),\quad m\in[-1,1],
\end{equation}
where $a=b'(0)=2\beta-2$ and $|Q(m)|\le K m^2$ for  $m\in[-1,1]$, we can write
\[
M_N(t)= a \int_0^t M_N(s)ds + \int_0^tQ(M_N(s))ds + \int_0^tdZ_N(s),
\]
We can now use variation of constants to write
\begin{equation}
\label{eq:var_of_const}
M_N(t)=e^{at}\int_0^t e^{-as}dZ_N + e^{at}\int_0^t e^{-as}Q(M_N(s))ds,
\end{equation}
where the integral w.r.t.\ $Z_N$ is understood as Lebesgue--Stieltjes integral.

\begin{lemma}\label{lm:clt} Suppose there is a sequence of stopping times $\theta_N$
satisfying
\begin{equation}
\label{eq:theta_to_infty}
\theta_N\stackrel{\Pp}{\to}\infty,\quad N\to\infty,
\end{equation}
and
\begin{equation}
\sup_{t\le \theta_N}|M_N(t)|\le N^{-\gamma},
\label{eq:M_stays_close_to_0}
\end{equation}
for some $\gamma>0$ and all $N$.
Then,  as $N\to\infty$,  $I_N=N^{1/2}\int_0^{\theta_N}e^{-as}dZ_N(s)$ converges in distribution to $\Nc(0,2a^{-1})$.
\end{lemma}

\bpf
Let us introduce a process
\begin{equation}
V_N(t)=\int_0^{t\wedge \theta_N}e^{-as}dZ_N(s).  
\label{eq:integral_martingale}
\end{equation}
It is a martingale with 
\[
[V_N]_t=\int_0^{t\wedge \theta_N}e^{-2as}d[Z_N]_s, 
\]
where square brackets denote the quadratic variation process.
Let us define $U_N(s)=V_N(f(s))$, where 
\[
f(s)=-\frac{\ln(1-as/2)}{2a},\quad s\in[0,2a^{-1}]. 
\]
We need the following statement which is a specific case of Theorem 1.4 in \cite[Chapter 7]{Ethier-Kurtz:MR838085}
(see also bibliographical notes therein for the history of this theorem and related results):
\begin{theorem}\label{th:invariance_principle}
 For each $N\in\N$, let $U_N$ be a martingale w.r.t.\ some filtration, with cadlag paths and $U_N(0)=0$. Suppose
 for all $t\in[0,2a^{-1}]$, $A_N(t)=[U_N]_t$ satisfies
\[
 N A_N(t)\stackrel{\Pp}{\to} t,\quad N\to\infty,
\]
and 
\begin{equation}
\label{eq:condition_on_smallness_of_jumps}
\lim_{N\to\infty} \E \left[N^{1/2}\sup_{0\le t\le 2a^{-1}}|U_N(t)-U_N(t-)|\right]=0.  
\end{equation}
Then, as $N\to\infty$, $N^{1/2}U_N$ converges in distribution in the Skorokhod topology to the standard Wiener process on $[0,2a^{-1}]$.
\end{theorem}

In our case, condition~\eqref{eq:condition_on_smallness_of_jumps} is fulfilled automatically since all the jumps of~$U_N$ are bounded by $2N^{-1}$ in absolute value. All the jumps
of $Z_N$ are equal to $2N^{-1}$ in absolute value, so that
\[
 [Z_N](t)=\frac{4}{N^2}B_N(t),\quad t\ge 0,
\]
where $B_N(t)$ denotes the number of jumps the process $Z_N$ makes up to time~$t$.
Next,
\[
 A_N(t)=\frac{4}{N^2}\sum_{\substack{s: s\le f(t)\wedge \theta_N\\ M_N(s)\ne M_N(s-)}} e^{-2as},
\]
and
\begin{equation}
 N A_N(t)\stackrel{\Pp}{\to} 4 \int_0^{f(t)} e^{-2as}ds= t.
\label{eq:convergence_of_A} 
\end{equation}
This is a consequence of the following result:
\begin{lemma} For any non-increasing function $g:\R_+\to\R_+$ and any $t\ge 0$,
\[
\frac{1}{N}\sum_{\substack{s: s\le t\wedge \theta_N\\ M_N(s)\ne M_N(s-)}} g(s)
\stackrel{\Pp}{\to} \int_0^{t} g(s)ds.
\]
\end{lemma}
\bpf It is sufficient to  assume that
$g(s)=\ONE_{s\in[0,h]}$ for some $h>0$ since one can use linear combinations of functions of this
form to approximate any non-increasing function. We have to show that for all $t$,
\[
 \frac{1}{N}B_N(t\wedge h\wedge \theta_N)\stackrel{\Pp}{\to} t\wedge h.
\]
Due to~\eqref{eq:theta_to_infty}, we can restrict ourselves to the high probability event $\{t<\theta_N\}$.
Let $(s_{N,i})$ be the increasing sequence of times of jumps for each $N$, i.e., 
$M_N(s_{N,i})\ne M_N(s_{N,i}-)$. Clearly, conditioned on $M_N(s_{N,i})$, the spacing
random variable $s_{N,i+1}-s_{N,i}$ has exponential distribution with parameter 
\[
\lambda_+(M_N(s_{N,i}),N)+\lambda_N(M_N(s_{N,i}),N)=\lambda(M_N(s_{N,i}),N),
\]
where
\[
\lambda(m,N)=\lambda_-(m,N)+\lambda_+(m,N).
\]
Therefore, using \eqref{eq:M_stays_close_to_0} we see that
the number of points $s_i$ not exceeding $t\wedge h$ is between two Poisson processes
with intensities $\lambda(0,N)=N$ and $\lambda(N^{-\beta},N)$ evaluated at time $t\wedge h$.
Our claim follows since $\lambda(N^{-\beta},N)/N\to1$.
\epf

Now all the conditions of Theorem~\ref{th:invariance_principle} have been verified and we conclude that $N^{1/2}U_N$ converges
in distribution in Skorokhod topology to the Wiener process on $[0,2a^{-1}]$. Therefore,
$ I_N=N^{1/2}U_N(2a^{-1})$ converges in distribution to $\Nc(0,2a^{-1})$, and the proof of Lemma~\ref{lm:clt} is complete. \epf

Let us now check that the conditions of Lemma~\ref{lm:clt} hold for
$ \theta_N =\inf\{t: |M_N(t)|\ge N^{-\gamma}\}$.
Notice that condition~\eqref{eq:M_stays_close_to_0}
is satisfied automatically.

\begin{lemma} If $\frac{1}{4}<\gamma<\frac{1}{2}$, then \eqref{eq:theta_to_infty} holds.
\end{lemma} 
\bpf Let us take any $T>0$.
On $\{\theta_N <T\}$ 
\begin{align}
 N^{-\gamma}&\le |M_N(\theta_N)|\notag \\
 &=e^{a\theta_N} \left| V_N(\theta_N)+\int_0^{\theta_N}e^{-as}Q(M_N(s))ds\right|\label{eq:equation_on_exit_gamma}\\
 &\le e^{aT} (|V_N(\theta_N)|+a^{-1}KN^{-2\gamma}),\notag
\end{align}
where $V$ has been introduced in~\eqref{eq:integral_martingale}.
The quadratic variation of the martingale $V_N$ was computed in the proof of Lemma~\ref{lm:clt}, 
and we can conclude that
\[
 \E [V_N]_{\theta_N}\le 2a^{-1} N^{-1}.
\]
Therefore, by the Chebyshev inequality and the fact that $V^2_N-[V_N]$ is a martingale (see 
Proposition~6.1 in \cite[Chapter 2]{Ethier-Kurtz:MR838085}), 
\begin{align*}
 \Pp\{|V_N(\theta_N)|>N^{-\gamma/2-1/4}\} &\le \frac{\E V^2_N(\theta_N)}{N^{-\gamma-1/2}}
\le \frac{\E [V_N]_{\theta_N}}{N^{-\gamma-1/2}}\le \frac{2 a^{-1}N^{-1}}{N^{-\gamma-1/2}} \to 0.
\end{align*}
Now, on $\{\theta_N <T\}\cap \{|V_N(\theta_N)|\le N^{-\gamma/2-1/4}\}$, we have
\[
 N^{-\gamma}\le e^{aT}( N^{-\gamma/2-1/4} + a^{-1}KN^{-2\gamma})
\]
which is impossible for large $N$ under our assumptions. We conclude that $P\{\tau_N<T\}\to 0$, and the lemma follows.
\epf
\medskip

We are ready to describe the asymptotics of the exit from $[-N^{-\gamma},N^{-\gamma}]$.
\begin{lemma} \label{lm:exit_diminishing_ball} If $1/4<\gamma<1/2$, then
\[
 \left(\sgn M_N(\theta_N),\ \theta_N -\frac{1/2-\gamma}{a}\ln N\right)\ \stackrel{distr}{\longrightarrow}\ \left(\sgn H, -\frac{1}{a}\ln|H|\right),
\]
where $\Law(H)= \Nc(0,2a^{-1})$.
\end{lemma}

\bpf Considering the process $M_N$ at time $\theta_N$ and using~\eqref{eq:var_of_const}, we obtain
\begin{align*}
 \theta_N&=\frac{1}{a}\ln\frac{\left\lceil N^{-\gamma}\cdot \frac{N}{2}\right\rceil\cdot\frac{2}{N}}{\left| V_N(\theta_N)+\int_0^{\theta_N}e^{-as}Q(M_N(s))ds\right|}\\
&=\frac{1/2-\gamma}{a}\ln N -\frac{1}{a}\ln\left|I_N+N^{1/2}\int_0^{\theta_N}e^{-as}Q(M_N(s))ds\right|+o(1).
\end{align*}
Also,
\[
\sgn M_N(\theta_N)=\sgn \left(I_N+N^{1/2}\int_0^{\theta_N}e^{-as}Q(M_N(s))ds\right).
\]
Since
\[
 N^{1/2}\int_0^{\theta_N}e^{-as}Q(M_N(s))ds\le \frac{N^{1/2}K}{a}N^{-2\gamma}\to 0,
\]
the desired statement follows from Lemma~\ref{lm:clt}.
\epf

\medskip

 Let us now study the exit of $M_N$ from an interval $[-r,r]$ where $r$ is a small number that does not depend on $N$.
We define $Y_N(t)=M_N(t+\theta_N)$ and for any $r\in(0,R)$ we define 
\[
\nu_N(r)=\inf\{t\ge 0:\ |Y_N(t)|= r\}. 
\]
We are going to compare the evolution of the magnetization process to the deterministic trajectory of the flow 
$S^t$ generated by the drift $b$.

For any $\delta>0$ we introduce  $t(\delta,r)$ as the only solution $t$ of $S^{t}\delta=r$, i.e., it is the time it takes
for the solution of ODE $\dot x=b(x)$ to travel from $\delta$ to $r$.

\begin{lemma}\label{lm:asymptotics_for_deterministic_exit_time}
For any $r>0$,
\[
\lim_{\delta\to 0} \left(t(\delta,r)-\frac{1}{a}\ln\frac{r}{\delta}\right)=K(r),
\]
where $K(\cdot)$ was defined in~\eqref{eq:K(R)}.
\end{lemma}
\bpf
By the basic formula for solutions of autonomous ODE's (see, e.g., \cite[Section 1.2]{Arnold:MR2242407}):
\begin{align*}
 t(\delta,r)&=\int_{\delta}^r\frac{dx}{b(x)}
=\int_{\delta}^r\left(\frac{1}{b(x)}-\frac{1}{ax}\right)dx + \int _{\delta}^r\frac{dx}{ax}\\
&=\int_{\delta}^r\frac{1}{b(x)}\frac{ax-b(x)}{ax}dx + \frac{1}{a}(\ln r-\ln\delta),
\end{align*}
and the lemma follows.\epf

\begin{lemma}\label{lm:follows_deterministic_in_small_fixed_ball}
There is  $r_0>0$ such that
\[
 \sup_{0\le t\le t(|Y_N(0)|,r_0)} |Y_N(t)-S^t(Y_N(0))|\stackrel{\Pp}{\to} 0,\quad N\to\infty.
\]
\end{lemma}
\bpf
Denote $\Delta_N(t)=Y_N(t)-S^tY_N(0)$. Since
\[
Y_N(t)=Y_N(0)+\int_0^tb(Y_N(s))ds+Z'_N(t),
\]
where $Z'_N(t)=Z_N(t+\theta_N)-Z_N(\theta_N)$ is a martingale, and
\[
S^tY_N(0)=Y_N(0)+ \int_0^tb(S^sY_N(0))ds,
\]
we see that for any $r\in (0,m_*)$ and any $t\in(0,t(|Y_N(0)|,r))$,
\[
|\Delta_N(t\wedge \nu_N(r))|\le L(r)\int_0^{t\wedge\nu_N(r) }|\Delta_N(s)|ds+\sup_{s\le t\wedge \nu_N(r)}|Z'_N(s)|, 
\]
where $L(r)$ is the Lipschitz constant of $b$ on $[-r,r]$.

Since $[Z'_N]_t=4N^{-2}(B_N(\theta_N+t)-B_N(\theta_N))$, and the number of jumps
between $\theta_N$ and $\theta_N+t$ is stochastically dominated by the increment of the
Poisson process with intensity $N$, we have
\[
 \Pp\left\{\sup_{s\le t(|Y_N(0)|,r)\wedge \nu_N(r)}|Z'_N(s)|> N^{-\delta}\right\}\to 0
\]
if $\delta<1/2$. On the complementary event, applying Gronwall's inequality, we obtain
for some constant $C>0$,
\begin{equation}
\label{eq:Delta_estimate}
|\Delta_N(t\wedge \nu_n(r))|\le e^{L(r) t(|Y_N(0)|,r)} N^{-\delta}\le C N^{\gamma L(r)/a}N^{-\delta}.
\end{equation}
We can choose $r$ to be so small that $L(r)$ is close to $a$ enough to ensure that $\gamma L(r)/a <1/2$.
Consequently, we can choose $\delta<1/2$ such that $\rho=\delta-\gamma L(r)/a>0$, and the r.h.s.\ of
\eqref{eq:Delta_estimate} converges to 0. For any $r_0\in(0,r)$ we conclude then that
$\Pp\{\nu_N(r)< t(|Y_N(0)|,r_0)\}\to 0$, and
\[
 \Pp\left\{\nu_N(r)\ge t(|Y_N(0)|,r_0);\  \sup_{s\le t(|Y_N(0)|,r_0)} |\Delta_N(s)|>N^{-\rho}\right\}\to 0.
\]
which completes the proof of the lemma.
\epf

\medskip

We can now combine the results of Lemmas~\ref{lm:exit_diminishing_ball} and~\ref{lm:follows_deterministic_in_small_fixed_ball}.

\begin{lemma} \label{lm:asymptotics_for_exit_from_fixed_small_ball}
For any $r\in(0,r_0)$,
\[
 \left(\sgn M_N(\tau_N(r)),\ \tau_N(r) -\frac{1}{2a}\ln N\right)\ \stackrel{distr}{\longrightarrow}\ \left (\sgn H, -\frac{1}{a}\ln|H|+\frac{\ln r}{a}+K(r)\right),
\]
\end{lemma}
\bpf Obviously, 
\begin{equation}\label{eq:tau_is_sum_of_stopping_times} 
\tau_N(r)=\theta_N+\nu_N(r). 
\end{equation}
Lemma~\ref{lm:follows_deterministic_in_small_fixed_ball} implies that
\[
 \nu_N(r)-t(|Y_N(0)|,r)\stackrel{\Pp}{\to} 0.
\]
This together with Lemma~\ref{lm:asymptotics_for_deterministic_exit_time} implies
\begin{equation}
\label{eq:asymptotics_for_nu}
 \nu_N(r)-\frac{1}{a}\ln\frac{r}{N^{-\gamma}}-K(r) \stackrel{\Pp}{\to} 0. 
\end{equation}
The lemma follows now from~\eqref{eq:tau_is_sum_of_stopping_times},\eqref{eq:asymptotics_for_nu},  and Lemma~\ref{lm:exit_diminishing_ball}.
\epf

The next result follows from the same considerations as Lemma~\ref{lm:follows_deterministic_in_small_fixed_ball}, except that it is easier since we consider
a finite time horizon.
\begin{lemma}\label{lm:estimate_on_finoite_time_horizon} Let $r$ be as in the last lemma.
 Let $Y_N(t)=M_N(\tau_N(r)+t)$. Then, for any $T>0$,
\[
 \sup_{0\le t\le T} |Y_N(t)-S^t(Y_N(0))|\stackrel{\Pp}{\to} 0,\quad N\to\infty.
\] 
\end{lemma}

This lemma means that after $\tau_N(r)$ the process essentially follows the deterministic trajectory. Since 
$H=2^{1/2}a^{-1/2} G$, where $G$ is standard Gaussian,
our main result is a
direct consequence of Lemmas~\ref{lm:asymptotics_for_exit_from_fixed_small_ball} and~\ref{lm:estimate_on_finoite_time_horizon}. In fact, it extends 
Lemma~\ref{lm:asymptotics_for_exit_from_fixed_small_ball} since the latter is valid only for sufficiently small values of threshold, whereas
our main result applies to any $R\in(0,m_*)$.
\epf

{\bf Acknowledgment.} The author is grateful to National Science Foundation for partial support through CAREER grant DMS-0742424.

\bibliographystyle{plain}
\bibliography{happydle}

\end{document}